
\tolerance=10000
\magnification=1200
\raggedbottom

\baselineskip=15pt
\parskip=1\jot

\def\sk{\vskip 3\jot}

\def\heading#1{\vskip3\jot{\noindent\bf #1}}
\def\label#1{{\noindent\it #1}}
\def\QED{\hbox{\rlap{$\sqcap$}$\sqcup$}}


\def\ref#1;#2;#3;#4;#5.{\item{[#1]} #2,#3,{\it #4},#5.}
\def\refinbook#1;#2;#3;#4;#5;#6.{\item{[#1]} #2, #3, #4, {\it #5},#6.} 
\def\refbook#1;#2;#3;#4.{\item{[#1]} #2,{\it #3},#4.}


\def\({\bigl(}
\def\){\bigr)}


\def\th{\vartheta}

\def\rh{\varrho}

\def\ph{\phi}

\def\ps{\psi}

\def\Th{\Theta}

\def\Ph{\Phi}
\def\Ps{\Psi}



\def\bfR{{\bf R}}

\input xy
\xyoption{all}

\def\abs#1{\big\vert#1\big\vert}
\def\bigabs#1{\left\vert#1\right\vert}

\def\expa{\exp\left(-{r^2\over n}\right)}
\def\expb{\exp\left(-{nx^2\over 4}\right)}
\def\expc{\exp\left(-{r^2\over 4a}\right)}

\def\Ex{{\rm Ex}}
\def\Var{{\rm Var}}
\def\Covar{{\rm Covar}}

\def\cotan{{\rm cotan}}
\def\cosec{{\rm cosec}}

{
\pageno=0
\nopagenumbers
\rightline{\tt 2d.2013.tex}
\vskip1in

\centerline{\bf Self-intersections of Two-Dimensional Equilateral Random Walks and Polygons}
\vskip0.5in

\centerline{Max B. Kutler}
\centerline{\tt kutler@uoregon.edu}
\sk

\centerline{Department of Mathematics}
\centerline{University of Oregon}
\centerline{Eugene, OR 97403}
\sk

\centerline{Margaret Rogers}
\centerline{\tt mrogers@hmc.edu}
\sk

\centerline{Nicholas Pippenger}
\centerline{\tt njp@math.hmc.edu}
\sk

\centerline{Department of Mathematics}
\centerline{Harvey Mudd College}
\centerline{301 Platt Bouldevard}
\centerline{Claremont, CA 91711}
\vskip0.5in

\noindent{\bf Abstract:}
We study the mean and variance of the number of self-intersections of the equilateral isotropic random walk in the plane, as well as the corresponding quantities for isotropic equilateral random polygons (random walks conditioned to return to their starting point after a given number of steps).
The expected number of self-intersections is $(2/\pi^2)n\log n + O(n)$ for both walks and polygons with $n$ steps.
The variance is $O(n^2 \log n)$ for both walks
and polygons, which shows that the number of self-intersections exhibits concentration around the mean.
\vfill\eject
}

\heading{1. Introduction}

The main objects of study in this paper are random walks and polygons in two dimensions.
We shall, however, also need to refer to  projections onto two dimensions of random walks and polygons in higher-dimensional spaces, so we shall begin by defining these objects in $d$-dimensional space.
A sequence $0 = X_0, X_1, \ldots, X_n$ of points in $\bfR^d$ will be called a {\it random walk\/} in $d$-dimensional space if the differences $X_1 - X_0, X_2 - X_1, \ldots, X_n - X_{n-1}$ are independent identically distributed random variables in $\bfR^d$.
We may also refer to the union $[X_0,X_1] \cup [X_1,X_2]\cup\cdots[X_{n-1},X_n]$ of the line segments $[X_{k-1},X_k]$ between the successive points $X_{k-1}$ and $X_k$ as the random walk.
All of the random walks we study will be {\t isotropic}; that is, the directions $(X_k - X_{k-1})/\abs{X_k - X_{k-1}}$ of the steps will always be uniformly distributed over the $(d-1)$-dimensional unit sphere, independent of the 
lengths $\abs{X_k - X_{k-1}}$ of the steps, so that the distribution of a random walk of length $n$ in $\bfR^d$ can be specified by giving the common distribution of the lengths.
A {\it random polygon\/} in $\bfR^d$ is a random walk in $\bfR^d$ conditioned on the event
$X_n = X_0$ of returning to the origin after $n$ steps.

A random walk or polygon in $\bfR^d$ with $d\ge 3$ can be projected onto a plane to give a random walk or polygon in $\bfR^2$.
Since our random walks are isotropic, the distribution of a projected walk or polygon will not depend on the choice of the plane onto which it is projected.

We shall be interested in the distribution of the number of {\it self-intersections\/} of a random walk or polygon.
Because angles are continuously distributed in our models, we can ignore the possibility that two points in a random walk coincide, that a point falls on a line segment, or that two line segments overlap in an interval of strictly positive length, since these events occur with probability zero.
When $n\ge 3$, the same observation applies to random polygons.
Thus the self-intersections occur at the interiors of line segments, and the number of 
self-intersections is the number of pairs of distinct line segments $[X_{k-1},X_k]$ and 
$[X_{l-1},X_l]$, with $1\le k<l\le n$ and $l-k\ge 2$.
that intersect at an interior point of each segment.
(In the case of polygons, we also exclude $l-k = n-1$.)

Diao and Ernst [D2] have studied the number of self-intersections of Gaussian random walks and polygons, showing that its mean is $(1/2\pi)n\log n + O(n)$ for both walks and polygons.
(For a Gaussian random walk, each step has an isotropic multivariate Gaussian distribution.
The number of self-intersections does not depend on the variance of the steps, since
this only affects the walk or polygon by a scale factor.
And, since the projection of a Gaussian random walk onto a smaller number of dimensions is again a Gaussian random walk, their result does not depend on the dimension of the original walk.)

Diao {\it et al.} [D1] have studied the corresponding problem for the projections onto two dimensions of three-dimensional equilateral random walks and polygons, obtaining
the estimate $(3/16)n\log n + O(n)$ for both walks and polygons.
(For an equilateral random walk, each step has unit length.
The projection of a three-dimensional equilateral walk onto two dimensions is not equilateral,
so their analysis is done in three-dimensional space.)

In this paper we study two-dimensional equilateral isotropic random walks and polygons.
As might be expected, our result for the mean number of self-intersections differs from the 
results cited above only in the constant factor in the leading term: we show that it is 
$(2/\pi^2)n\log n + O(n)$ for both walks and polygons.
But we carry the analysis further than that of the results cited above, and show that the variance
is $O(n^2 \log n)$ for both walks and polygons.
Thus the number of self-intersections exhibits concentration about its mean in both cases.
Indeed, by Chebyshev's inequality, the probability that the number of self-intersections
differs from its mean by more than $n (\log n)^{1/3}$ is at most $O\(1/(\log n)^{1/3}\)$.
Finally, we observe that $1/2\pi = 0.1591\ldots\,$, $3/16 = 0.1875$, and $2/\pi^2 = 0.2026\ldots\,$.
Thus equilateral isotropic random walks and polygons have on the average more self-intersections than their counterparts in either of the other models mentioned above.
\sk

\heading{2. Quasi-Gaussian Densities}

We shall say that a two-dimensional probability density $f_{R,\Th}(r,\th)$
is {\it $n$-Gaussian\/} if
its polar coordinates $(R,\Th)$ have a density of the form
$$f_{R,\Th}(r,\th) = {1\over \pi n} \expa.$$
In this paper we shall often encounter two-dimensional densities that are approximately, but not exactly, Gaussian.
In this section we shall define a suitable notion of ``approximately Gaussian'' which we call 
``quasi-Gaussian''.
We shall show that the sum of $n$ steps of an equilateral isotropic random walk
is quasi-Gaussian.

We shall say that a two-dimensional probability density $f_{R,\Th}(r,\th)$
is {\it $n$-quasi-Gaussian\/} if
its polar coordinates $(R,\Th)$ have a density of the form
$$f_{R,\Th}(r,\th) = {1\over \pi n} \expa  + O\left({1\over n^2}\right), \eqno(2.1)$$
where the constant in the $O$-term is independent of $r$, $\th$ and $n$.
(This definition really applies to a family of densities parameterized by $n$, and
it only refers to their behavior for large $n$.
It does not require that $f_{R,\Th}(r,\th)$ be independent of $\th$,
but only that its dependence on $\th$  affects the density by at most $O(1/n^2)$.)

The distribution of the isotropic equilateral random walk was apparently first treated by Rayleigh [R1, pp.~35--42] in 1877.
This random walk gives the distribution of the amplitude and phase of the sum of identical sinusoidal oscillations with equal amplitudes and random phases.
Rayleigh gave the asymptotic formula
$$f^*_R(r) \sim {2r\over n} \, \expa,$$
for the radial density function of this walk after $n$ steps.
Since the walk is isotropic, dividing by $2\pi r$ gives the density
$$f^*_{R,\Th}(r,\th) \sim {1\over \pi n} \expa,$$
which agrees with the $n$-Gaussian factor in (2.1).

In 1906, Kluyver [K] gave the integral representation
$$F^*_R(r) = r \int_0^\infty J_1(rx) \, J_0(x)^n \, dx \eqno(2.2)$$
for the radial distribution function, where $J_n(x)$ is the Bessel function of order $n$
(see Watson [W1]).
The representation
$$f^*_R(r) = r \int_0^\infty J_0(rx) \, J_0(x)^n \, x \, dx, \eqno(2.3)$$
for the radial density function can be obtained from (2.2) by differentiating with respect to $r$,
then using the identity $J'_0(x) = -J_1(x)$ (see Watson [W1], p.~18)
and the differential equation $x J''_0(x) + J'_0(x) + x J_0(x)  =0$ (see Watson [W1], p.~19):
$$\eqalign{
{d\over dr} r J_1(rx)
&= -{d\over dr} r J'_0(rx) \cr
&= - rx J''_0(rx) - J'_0(rx) \cr
&= rx J_0(rx). \cr
}$$
(For $n\ge 5$ the integral is absolutely convergent: we have $\abs{J_1(rx)}\le 1$ (see Watson [W1, p.~31]) and $J_0(x) = O(1/x^{1/2})$ (see Watson [W1, p.~195]).
This fact justifies the differentiation of the integral (see for example Whittaker and Watson [W2, 
p.~174]).)
Dividing (2.3) by $2\pi r$ yields
$$f^*_{R,\Th}(r,\th) = {1\over 2\pi} \int_0^\infty J_0(rx) \, J_0(x)^n \, x \, dx, \eqno(2.4)$$

In 1919, Rayleigh [R2] gave a heuristic derivation of an asymptotic expansion for $f^*_R(r)$ that, if proved rigorously, would show that the equilateral random walk is quasi-Gaussian.
(Rayleigh's derivation involves differentiating an asymptotic expansion term-by-term.)
We shall give a rigorous proof below that the equilateral random walk is $n$-quasi-Gaussian.
(Our proof could be extended to establish the complete asymptotic expansion, up to terms of order $O(1/n^k)$ for any fixed $k\ge 1$.)

\label{Proposition 2.1}
The sum of $n$ steps of an equilateral isotropic random walk is $n$-quasi-Gaussian.

\label{Proof:}
In view of (2.4), it will suffice to show that
$${1\over 2\pi} \int_0^\infty J_0(rx) \, J_0(x)^n \, x \, dx
= {1\over \pi n} \expa + O\left({1\over n^2}\right), \eqno(2.5)$$
where the constant in the $O$-term is independent of both $r$ and $n$.
Let $x_0 = (24\log n / n)^{1/2}$.
Our first step will be to show that 
$$ \int_{x_0}^\infty J_0(rx) \, J_0(x)^n \, x \, dx
=  O\left({1\over n^2}\right). \eqno(2.6)$$
Let $x_1 = n^{1/2}$.
We shall prove (2.6) by showing that
$$ \int_{x_0}^{x_1} J_0(rx) \, J_0(x)^n \, x \, dx
=  O\left({1\over n^2}\right) \eqno(2.7)$$
and
$$ \int_{x_1}^\infty J_0(rx) \, J_0(x)^n \, x \, dx
=  O\left({1\over n^2}\right). \eqno(2.8)$$
To prove (2.7), we first observe that $J_0(x)$ is analytic for $x\in[0,\infty)$ and 
$J_0(x) = 1 - x^2/4 + O(x^4)$  as $x\to 0$ (see Watson [W1, p.~16]), and that $J_0(x)$ assumes values near $1$ only for $x$ near $0$.
(The last fact can easily be seen from a graph of $J_0(x)$; we shall indicate how it can be derived from facts  proved by Watson [W1], who gives no graphs!)
Firstly, the integral representation
$$J_0(x) = {1\over\pi} \int_0^\pi \cos (x\cos \th)\,d\th$$
(see Watson [W1, p.~24]) shows that $J_0(x) = 1$ only for $x=0$ (because it is an average of quantities that are all $1$ only if $x=0$).
Thus we cannot have $J_0(x_n)\to 1$ for $x_n\to x$ for any finite $x>0$, because $J_0(x)$, being analytic, is continuous.
We also cannot have $J_0(x_n)\to 1$ for $x_n\to \infty$, since $J_0(x)\to 0$ as $x\to\infty$
(see Watson [W1, p.~195]).
Thus $J_0(x)$ assumes values near $1$ only for $x$ near $0$.)
Thus there exists $\xi_0>0$ such that, for $x\le \xi_0$, we have not only  $\abs{J_0(x)} \le 1 - x^2/8$, 
but also $\abs{J_0(y)} \le 1 - x^2/8$ for all $y\ge x$.
We have $x_0\le \xi_0$ for all sufficiently large $n$.
Since $\abs{J_0(rx)}\le 1$, we then have
$$\eqalign{
 \bigabs{\int_{x_0}^{x_1} J_0(rx) \, J_0(x)^n \, x \, dx}
 &\le x_1^2 \, (1 - x_0^2 / 8)^n \cr
 &\le x_1^2 \, \exp\left(-{nx_0^2\over 8}\right) \cr
 &=  O\left({1\over n^2}\right), \cr
}$$
which proves (2.7).
To prove (2.8), we observe that since $J_0(x) = (2/\pi x)^{1/2} \, \cos (x-\pi/4) + O(1/x)$
(see Watson [W1, p.~195]), there exists
$\xi_1$ such that $\abs{J_0(x)}\le 1/x^{1/2}$ for all $x\ge \xi_1$.
We have $x_1\ge\xi_1$ for all sufficiently large $n$.
If in addition we have $n\ge 16$, we then have
$$\eqalign{
 \bigabs{\int_{x_1}^\infty J_0(rx) \, J_0(x)^n \, x \, dx}
 &\le \int_{x_1}^\infty {dx\over x^{n/2-1}} \cr
&=  O\left({1\over x_1^{n/2-2}}\right) \cr
&=  O\left({1\over x_1^6}\right) \cr
&=  O\left({\log^3 n\over n^3}\right), \cr
 }$$
 which proves (2.8), and completes the proof of (2.6).
 Thus to prove (2.5), it will suffice to show that
 $${1\over 2\pi} \int_0^{x_0} J_0(rx) \, J_0(x)^n \, x \, dx
= {1\over \pi n} \expa + O\left({1\over n^2}\right). \eqno(2.9)$$
 
 Our next step will be to estimate the factor $J_0(x)^n$ of the integrand in (2.9) over the range 
 $x\in[0,x_0]$.
 Since $J_0(x) = 1 - x^2/4 + x^4/64 + O(x^6)$ in this range, $\log (1+y) = y - y^2/2 + O(y^3)$ for $x<1$, and
 $\exp z = 1 + z + z^2/2 + O(z^3)$ as $z\to 0$, we have
 $$\eqalign{
 J_0(x)^n
 &= \exp \(n\log J_0(x)\) \cr
 &= \exp \left(n\log \left( 1 - {x^2\over 4} + {x^4\over 64} + O(x^6)\right)\right) \cr
  &= \exp \left( -{nx^2\over 4} - {nx^4\over 64} + O(nx^6)\right) \cr
  &= \expb \, \exp \left(  - {nx^4\over 64} + O(nx^6)\right) \cr
  &= \expb \, \left(1 - {nx^4\over 64} + O(nx^6) + O(n^2 x^8)\right). \cr
 }$$
 Thus we have
 $$\eqalignno{ \int_0^{x_0} J_0(rx) \, J_0(x)^n \, x \, dx 
&=  \int_0^{x_0} J_0(rx) \,  \expb  \, x \, dx \cr
&\qquad - {n\over 64}\int_0^{x_0} J_0(rx) \,  \expb \,  x^5 \, dx \cr
&\qquad + O\left(\int_0^{x_0} nx^7 \, dx\right)
+ O\left(\int_0^{x_0} n^2 x^9 \, dx\right). &(2.10)\cr
 }$$
Since $nx_0^8 = O(\log^4 n/n^3)$ and $n^2 x^{10} = O(\log^5 n/n^3)$,
both of the last two terms in (2.10) are $O(1/n^2)$, so we have
$$\eqalignno{
\int_0^{x_0} J_0(rx) \, J_0(x)^n \, x \, dx 
&= \int_0^{x_0} J_0(rx) \,  \expb  \, x \, dx \cr
&\qquad - {n\over 64}\int_0^{x_0} J_0(rx) \,  \expb \,  x^5 \, dx
+ O\left({1\over n^2}\right). &(2.11) \cr
}$$
To simplify the integrals on the right-hand side of (2.11), we observe that
$$\eqalign{
\int_{x_0}^\infty  \expb \, x \, dx 
&= -{2\over n} \, \expb\bigg\vert_{x=x_0}^\infty \cr
&= O\left({1\over n^7}\right) \cr
}$$
and
$$\eqalign{
\int_{x_0}^\infty  \expb \, x^5 \, dx 
&= -\left({64\over n^3} + {16 x^2\over n^2} + {2x^4\over n}\right) \, \expb\bigg\vert_{x=x_0}^\infty \cr
&= O\left({\log^2 n\over n^7}\right) \cr
}$$
Since $\abs{J_0(x)}\le 1$, these estimates imply
$$\int_{x_0}^\infty  J_0(rx) \expb \, x \, dx  = O\left({1\over n^7}\right)$$
and
$$-{n\over 64}\int_{x_0}^\infty  J_0(rx) \expb \, x^5 \, dx = O\left({\log^2 n\over n^6}\right).$$
Adding these integrals to the right-hand side of (2.11), we obtain
$$\eqalignno{
\int_0^{x_0} J_0(rx) \, J_0(x)^n \, x \, dx 
&= \int_0^\infty J_0(rx) \,  \expb  \, x \, dx \cr
&\qquad - {n\over 64}\int_0^{x_0} J_0(rx) \,  \expb \,  x^5 \, dx
+ O\left({1\over n^2}\right). &(2.12) \cr
}$$
To evaluate the integrals on the right-hand side of (2.12), we shall use the integrals
$$\int_0^\infty J_0(rx) \, \exp(-ax^2) \, x \, dx = \expc \eqno(2.13)$$
(see Watson [W1, p.~393]) and
$$\int_0^\infty J_0(rx) \, \exp(-ax^2) \, x^5 \, dx 
= \left({1\over a^3} -{r^2\over 2a^4} + {r^4\over 32a^5}\right) \, \expc, \eqno(2.14)$$
which can be obtained from (2.13) by differentiating with respect to $a$
(see for example Whittaker and Watson [W2, p.~74]).
Applying these integrals to (2.12) with $a=n/4$ yields
$$\eqalignno{
\int_0^\infty J_0(rx) \, J_0(x)^n \, x \, dx 
&= {2\over n} \, \expa \cr
&\qquad - \left({1\over n^2} - {2r^2\over n^3} + {r^4\over 2n^4}\right) \, \expa
+ O\left({1\over n^2}\right). \cr
}$$
Since $r^2 \, \exp(-r^2/n) = O(n)$ and $r^4 \, \exp(-r^2/n) = O(n^2)$, we obtain
$$\eqalignno{
\int_0^\infty J_0(rx) \, J_0(x)^n \, x \, dx 
&= {2\over n} \, \expa
+ O\left({1\over n^2}\right). \cr
}$$
Multiplying by $1/2\pi$ yields (2.9), which completes the proof of (2.5).
\QED
\vfill\eject

\heading{3. A Triple Integral}

In this section we shall evaluate a triple integral that gives the coefficient $2/\pi^2$
of the $n\log n$ term in our results.
We consider the following geometric situation.
Let $Q$ denote the line segment of length $R$ from the origin to $(R,\Th)$,
where $(R,\Th)$ has an $m$-quasi-Gaussian density.
Let $S$ denote the line segment from the origin of unit length and making an angle $-\pi<\Ps<\pi$, measured 
counterclockwise from $Q$.
Let $T$ denote the line segment from $(R,\Th)$ of unit length making angle $-\pi<\Ph<\pi$, measured clockwise from $Q$.
Let $\Ps$ and $\Ph$ be uniformly distributed in the interval $(-\pi,\pi)$,
independently of each other and of $R$.
Let $E_m$ be the event that the segments $S$ and $T$ intersect at an interior point.

\label{Proposition 3.1:}
We have
$$\Pr[E_m] = {2\over \pi^2 m} + O\left({1\over m^2}\right). \eqno(3.1)$$

\label{Proof:}
Define the indicator function $I(r,\ps,\ph)$ to be $1$ (or $0$) according as the segments $S$ and $T$ do (or do not) intersect
at an interior point when $R$, $\Ps$ and $\Ph$ assume the values $r$, $\ps$ and $\ph$, respectively.
Then
$$\Pr[E_m] = {1\over (2\pi)^2}
\int_{-\pi}^\pi \int_{-\pi}^\pi \int_0^\infty I(r,\ps,\ph) \, f_{R}(r)\,dr\,d\ps\,d\ph.$$
Since $(R,\Th)$ is $m$-quasi-Gaussian, we have
$$\Pr[E_m] = {1\over (2\pi)^2}
\int_{-\pi}^\pi \int_{-\pi}^\pi \int_0^\infty I(r,\ps,\ph) \,
\left({2r\over m}\exp\left(-{r^2\over m}\right) + O\left({1\over m^2}\right)\right)\,dr\,d\ps\,d\ph.$$
Since $I(r,\ps,\ph)$ vanishes unless $r<2$, we can reduce the upper limit of the innermost integral:
$$\Pr[E_m] = {1\over (2\pi)^2}
\int_{-\pi}^\pi \int_{-\pi}^\pi \int_0^2 I(r,\ps,\ph) \,
\left({2r\over m}\exp\left(-{r^2\over m}\right) + O\left({1\over m^2}\right)\right)\,dr\,d\ps\,d\ph.$$
Since $\exp(x) = 1 + O(x)$ as $x\to 0$, we obtain
$$\eqalign{
\Pr[E_m] 
&= {1\over (2\pi)^2}
\int_{-\pi}^\pi \int_{-\pi}^\pi \int_0^2 I(r,\ps,\ph) \,
\left({2r\over m} + O\left({1\over m^2}\right)\right)\,dr\,d\ps\,d\ph \cr
&= {1\over 2\pi^2\,m}
\int_{-\pi}^\pi \int_{-\pi}^\pi \int_0^2 I(r,\ps,\ph) \,r\,dr\,d\ps\,d\ph 
+  O\left({1\over m^2}\right). \cr
}$$
Thus it will suffice to show that
$$\int_{-\pi}^\pi \int_{-\pi}^\pi \int_0^2 I(r,\ps,\ph)\,r\,dr\,d\ps\,d\ph = 4. \eqno(3.2)$$

Let $J$ denote the triple integral in (3.2). 
We shall show that $J = 4$.
It is clear that there is a function $\rh:(-\pi,\pi)\times(-\pi,\pi)\to [0,2]$ such that
$$I(r,\ps,\ph) = \cases{
1, &if $r\le \rh(\ps,\ph)$; \cr
0, &if $r > \rh(\ps,\ph)$. \cr
}$$
Thus
$$\eqalign{
J
&= \int_{-\pi}^\pi \int_{-\pi}^\pi \int_0^{\rh(\ps,\ph)} r\,dr\,d\ps\,d\ph \cr
& = {1\over 2}  \int_{-\pi}^\pi \int_{-\pi}^\pi \rh(\ps,\ph)^2\,d\ps\,d\ph. \cr
}$$ 
It is clear that $\rh(\ps,\ph)$ vanishes unless $\ps$ and $\ph$ have the same sign, and that it is unchanged if this common sign is reversed; thus
$$J =  \int_{0}^\pi \int_{0}^\pi \rh(\ps,\ph)^2\,d\ps\,d\ph.$$
It is clear that $\rh(\ps,\ph)$ vanishes if both $\ps$ and $\ph$ are obtuse 
(that is, belong to $(\pi/2,\pi)$).
Thus we may break the  integral into three parts, according as $\ps$, $\ph$, or neither is obtuse:
$$J =  \int_{\pi/2}^{\pi} \int_{0}^{\pi/2} \rh(\ps,\ph)^2\,d\ps\,d\ph + 
 \int_{0}^{\pi/2} \int_{\pi/2}^{\pi} \rh(\ps,\ph)^2\,d\ps\,d\ph + 
 \int_{0}^{\pi/2}  \int_{0}^{\pi/2} \rh(\ps,\ph)^2\,d\ps\,d\ph.$$
Since $\rh(\ps,\ph)$ is unchanged by the exchange of $\ps$ and $\ph$, the second term equals the first, so
$$J = 2 \int_{\pi/2}^{\pi} \int_{0}^{\pi/2} \rh(\ps,\ph)^2\,d\ps\,d\ph + 
 \int_{0}^{\pi/2}  \int_{0}^{\pi/2} \rh(\ps,\ph)^2\,d\ps\,d\ph.$$
 Furthermore, since $\rh(\ps,\ph)=0$ unless $\ps<\pi-\ph$, we can restrict the range of the inner integral in the first term, so
 $$J = 2 \int_{\pi/2}^{\pi} \int_{0}^{\pi-\ph} \rh(\ps,\ph)^2\,d\ps\,d\ph + 
 \int_{0}^{\pi/2}  \int_{0}^{\pi/2} \rh(\ps,\ph)^2\,d\ps\,d\ph.$$
 Again using the symmetry between $\ps$ and $\ph$, we may restrict the range of the inner integral in the second term to $\ps<\ph$, and double the resulting term, so
 $$
\eqalignno{ 
J 
&= 2 \int_{\pi/2}^{\pi} \int_{0}^{\pi-\ph} \rh(\ps,\ph)^2\,d\ps\,d\ph + 
2 \int_{0}^{\pi/2}  \int_{0}^{\ph} \rh(\ps,\ph)^2\,d\ps\,d\ph \cr
&= 2\int_0^{\pi/2} \int_\ps^{\pi-\ps} \rh(\ps,\ph)^2 \, d\ph \, d\ps. &(3.3) \cr
}$$

For the range of integration in (3.3), a little trigonometry shows that
$$\rh(\ps,\ph) = \sin(\ph+\ps) \, \cosec\,  \ph. \eqno(3.4)$$
To see this, we may imagine starting with $r=2$ and then reducing $r$ until $S$ and $T$ intersect,
which happens when $r=\rh(\ps,\ph)$.
Then in the resulting triangle, the side opposite $\ph$ has length $1$, while $\rh(\ps,\ph)$ is the length of the side opposite the angle $\pi-\ph-\ps$.
Applying the law of sines, and using the fact that $\sin(\pi-\ph-\ps) = \sin(\ph+\ps)$,   we obtain 
(3.4).
Substituting (3.4) into (3.3), we obtain
$$J = 2\int_0^{\pi/2} \int_\ps^{\pi-\ps} \sin^2 (\ph+\ps) \, \cosec^2 \ph \, d\ph \, d\ps. \eqno(3.5)$$
From the antiderivative
$$\int \sin^2 (\ph+\ps) \, \cosec^2 \ph \, d\ph
= \ph \cos^2 \ps - (\ph + \cotan\,\ph)\sin^2 \ps + \log \sin \ph \sin (2\ps),$$
we obtain
$$ \int_\ps^{\pi-\ps} \sin^2 (\ph+\ps) \, \cosec^2 \ph \, d\ph
= (\pi-2\ps)\cos(2\ps) + \sin(2\ps).$$
Substituting this value for the inner integral in (3.5) yields
$$J = 2\int_0^{\pi/2} (\pi-2\ps)\cos(2\ps) + \sin(2\ps) \, d\ps,$$
and evaluating this integral we obtain $J=4$ as desired.
This completes the proof of (3.2).
\QED

We observe for future reference that Proposition 3.1 continues to hold if a fixed constant
displacement $C$ is added to the $m$-quasi-Gaussian step from origin to $(R,\Th)$, because an $m$-quasi-Gaussian distribution assigns densities that differ at most by a factor $\(1 + O(1/m)\)$
to all points within a bounded distance of the origin.
(The constant in the $O$-term may now depend on $C$.)
\sk

\heading{4. The Mean}

We begin with walks.
Let the random variable $K_n$ denote the number of self-intersections in a two-dimensional equilateral isotropic random walk.
In this section, we shall derive the estimate
$$\Ex[K_n] = {2\over \pi^2} \; n\log n + O(n). \eqno(4.1)$$
Let $L_{i,j}$ denote the event that the $i$-th segment $[X_{i-1},X_i]$ intersects the $j$-th segment 
$[X_{j-1},X_j]$ at an interior point of each segment.
We shall also write $L_{i,j}$ for the indicator function of that event, assuming the value $1$ when the event occurs, and the value $0$ when it does not, so that $\Ex[L_{i,j}] = \Pr[L_{i,j}]$.
By the linearity of expectation, we have
$$\Ex[K_n] =  \sum_{1\le i<j\le n} \Pr[L_{i,j}]. \eqno(4.2)$$
Our problem is now to estimate $ \Pr[L_{i,j}]$.

It is clear that $ \Pr[L_{i,j}]$ depends only on the number $a = j-i-1$ of steps between the end of the $i$-th step and the beginning of the $j$-th step, and that it is zero unless $a\ge 1$.
Furthermore, $L_{i,j}$ has the same probability as the event $E_a$ defined in the preceding section
($(R,\Th)$ is the polar representation of the sum of the $(i+1)$-st through the $(j-1)$-st steps,
$S$ is the $i$-th step and $T$ is the $j$-th step).
Let $b = n-j-1$, as shown in Figure 4.1.

$$\xymatrix{
1 \ar@{-}[r] 
&i \ar@{-}[rr] _a \ar@{--}@(ur,ul)[rr] 
&&j \ar@{-}[r] _b
&n \cr
}$$
\centerline{Figure 4.1}
\sk

\noindent
By Proposition 3.1, we have
$$\Pr[L_{i,j}] = {2\over \pi^2\,a} + O\left({1\over a^2}\right).$$
Summing over the possible values of $a$ and $b$, we obtain
$$\eqalign{
\Ex[K_n] 
&= \sum_{0\le b\le n-3} \; \sum_{1\le a\le n-2-b} \Pr[L_{i,j}] \cr
&= \sum_{0\le b\le n-3} \; \sum_{1\le a\le n-2-b} \left({2\over \pi^2\,a} + O\left({1\over a^2}\right)\right) \cr
&= {2\over \pi^2}\,n\log n + O(n), \cr
}$$
because
$$\sum_{1\le m} {1\over v} = \log m + O(1)$$
and
$$\sum_{1\le v\le m} \log v = m\log m - m + O(\log m).$$
Thus (4.1) is verified.

We turn now to polygons.
Let the random variable $K'_n$ denote the number of self-intersections in a two-dimensional equilateral isotopic random polygon.
We shall derive the estimate
$$\Ex[K'_n] = {2\over \pi^2} \; n\log n + O(n). \eqno(4.3)$$

For polygons, with their circular symmetry, there is only one topological configuration for a self-intersection, as depicted in Figure 4.2.
\vskip0.25in
$$\xymatrix{
&\bullet \ar@{--}[rr] \ar@{-}@/_2.5pc/[rr] \ar@{-}@/^2.5pc/[rr]
&&\bullet \cr
}$$
\vskip0.35in
\centerline{\ \ \ \ \ \ \ \ \ \ Figure 4.2}
\sk

\noindent
The two ``stretches'' (indicated by solid lines) in Figure 4.2 are topologically equivalent,
and drawing wither of them at the bottom results in a picture like Figure 4.3.
\vskip0.15in
$$\xymatrix{
&i \ar@{-}[rr] _a \ar@{--}@(ur,ul)[rr] \ar@{-}@/_2.5pc/[rr] ^b
&&j  \cr
}$$
\vskip0.15in
\centerline{\ \ \ \ \ \ \ \ \ \ Figure 4.3}
\sk

\noindent
We shall analyze this case with arguments similar to those we used for Figure 4.1.
But for polygons we have the constraint $a+b = n-2$, so we shall not sum over
all combinations of $a$ and $b$.
Since at least one of the two stretches in Figure 4.2 must have length
at least $(n-2)/2\ge n/4$ (assuming $n\ge 4$), we may choose such a long stretch 
to draw at the bottom, so that we have $b\ge n/4$.
This choice implies $a\le (n-2)/2$.
Thus we shall sum $\Pr[L_{i,j}]$ over $1\le a\le (n-2)/2$ and then include an extra factor of $n$ to the sum, to take account of the $n$ possible 
positions in which this figure might appear around the polygon. 

By circular symmetry, $\Pr[L_{i,j}]$ depends only on $a = j-i-1$ and on $n$.
To determine this dependence, we must reconsider the situation described in the preceding section, conditioning on the event
that $T$ and $S$ are the first and last steps of an equilateral isotropic random walk of $b+2$ steps from $(R,\Th)$ back to the origin.
This introduces two complicating effects.
First, the density of $(R,\Th)$ is no longer $a$-quasi-Gaussian.
We shall see, however, that is it quasi-Gaussian with a smaller parameter.
Second, $\Ps$ and $\Ph$ are no longer independent and uniformly distributed.
We shall see, however, that they are close to being so.

We begin by reconsidering the density of $(R,\Th)$.
It is well known that the density of an $v$-Gaussian step, conditioned on the event that a further independent $w$-Gaussian step
returns to the origin, is $(v\parallel w)$-Gaussian, where $(v\parallel w) = vw/(v+w)$ is the ``harmonic sum'' of $v$ and $w$.
(The variances of ``parallel'' Gaussian steps combine like resistance in parallel.)
To derive this result, we have only to multiply the densities and integrate the result to renormalize.
For quasi-Gaussian steps, we must add error terms $O(1/v^2)$ and $O(1/w^2)$.
But since
$${1\over v^2} + {1\over w^2} = {1\over v^2\parallel w^2} = {v^2 + w^2 \over v^2 \, w^2}
\le {(v+w)^2\over (vw)^2} = {1\over (v\parallel w)^2},$$
these error terms can be combined into a single one of order $O\(1/(v\parallel w)^2\)$.
Since this is the error term for an $(v\parallel w)$-quasi-Gaussian density, we conclude that the density of $(R,\Th)$
for polygons is $\(a\parallel (b+2)\)$-quasi-Gaussian.

We observe for future reference that if fixed constant displacements $C$ and $D$ are added to the
$v$- and $w$-quasi-Gaussian steps considered above, their parallel connection is then
$(v\parallel w)$-quasi-Gaussian with an added constant displacement $(vC + wD)/(v+w)$.

We turn next to the dependence of $\Ps$ and $\Ph$.
Since $I(r,\ps,\ph)$ vanishes unless $r<2$, we assume that $R$ is some value $r<2$.
In this case the other endpoints of the unit segments $S$ and $T$ are at distance at most $4$.
Suppose in addition that $\Ps$ assumes some value $\ps$ (thus determining the position of the other endpoint $s$ of $S$).
An $(b+2)$-quasi-Gaussian density assigns to all points within distance at most $4$ of  $s$
values that differ at most by a factor $1 + O(1/b)$.
Thus we have $f_{\Ph\mid \Ps=\ps,R=r}(\ph) = (1/2\pi)\(1 + O\(1/b)\)$ when $r<2$.
By the same argument, we have $f_{\Ps\mid \Ph=\ph,R=r}(\ps) = (1/2\pi)\(1 + O\(1/b)\)$ when $r<2$.
Thus we have
$f_{\Ps,\Ph\mid R=r}(\ps,\ph) = (1/2\pi)^2\(1 + O\(1/b)\)$ when $r<2$.

Substituting this result for the factor of $1/(2\pi)^2$ that represents the uniform joint distribution of $\Ps$ and $\Ph$
in the preceding section, and changing the $a$-quasi-Gaussian distribution of $(R,\Th)$ to an 
$\(a\parallel (b+2)\)$-quasi-Gaussian distribution (as indicated in the preceding paragraph), we conclude that
$$\eqalign{
\Pr[L_{i,j}] 
&= \left({2\over \pi^2 (a\parallel (b+2))} + O\left({1\over (a\parallel (b+2))^2}\right)\right) \, 
\left(1 + O\left({1\over b}\right)\right) \cr
&= {2\over \pi^2 a} + O\left({1\over a^2}\right), \cr
}$$
where we have used the fact that $b\ge n/4$ to simplify the parallel combinations involving $b$,
together with $1/n = O(1/a)$.
Summing over $1\le a\le n/2$ and multiplying by an additional factor of $n$ then yields (4.3).
\sk

\heading{5. The Variance}

We again begin with walks.
Recall that the random variable $K_n$ denotes the number of self-intersections in a two-dimensional equilateral isotropic random walk.
We shall derive the estimate
$$\Var[K_n] = O(n^2 \log n) \eqno(5.1)$$
We shall use the formula
$$\Var[K_n] = \sum_{1\le i<j\le n \atop 1\le k<l\le n} \Covar[L_{i,j}, L_{k,l}], \eqno(5.2)$$
where
$$\Covar[L_{i,j},L_{k,l}] = \Pr[L_{i,j},L_{k,l}] - \Pr[L_{i,j}]\cdot\Pr[L_{k,l}].\eqno(5.3)$$

We shall suppose to begin with that $i$, $j$, $k$ and $l$ all distinct. 
In fact we shall suppose that they differ pairwise by at least $2$.
We consider first the terms with $i<j<k<l$, as depicted in Figure 5.1.
\vskip0.15in
$$\xymatrix{
1 \ar@{-}[r] 
&i \ar@{-}[rr] \ar@{--}@(ur,ul)[rr] 
&&j \ar@{-}[rr] 
&&k \ar@{-}[rr] \ar@{--}@(ur,ul)[rr] 
&&l \ar@{-}[r] 
&n \cr
}$$
\centerline{Figure 5.1}
\sk

\noindent 
For these terms, $L_{i,j}$ and $L_{k,l}$ are independent, so 
$\Pr[L_{i,j},L_{k,l}] = \Pr[L_{i,j}]\cdot\Pr[L_{k,l}]$ and
$\Covar[L_{i,j}, L_{k,l}]$ vanishes.
This observation of course also applies to terms with $k<l<i<j$.

We consider next the terms with $k<i<j<l$, as depicted in Figure 5.2.
\vfill\eject

\vskip0.55in
$$\xymatrix{
1 \ar@{-}[r] 
&k \ar@{-}[r] _a \ar@{--}@/^3.75pc/[rrrr]
&i \ar@{-}[rr] _b \ar@{--}@(ur,ul)[rr] 
&&j \ar@{-}[r] _c
&l \ar@{-}[r] _d
&n 
\cr
}$$
\centerline{Figure 5.2}
\sk

\noindent 
For these terms, we shall use the upper bound
$$\Covar[L_{i,j}, L_{k,l}] \le \Pr[L_{i,j},L_{k,l}] = \Pr[L_{i,j}] \, \Pr[L_{k,l} \mid L_{i,j}]. \eqno(5.4)$$
(Since $\Var[K_n]$ is non-negative, we can upper-bound it by upper-bounding
each term in the sum (5.2).)
We shall define $a = i-k-1$, $b = j-i-1$, $c = l-j-1$ and $d=n-l-1$.
Then we must sum over terms with $a, b, c, d\ge 1$ and $a+b+c+d\le n-4$.
By Proposition 3.1,
$$\Pr[L_{i,j}] = {2\over \pi^2 b} + O\left({1\over b^2}\right),$$
but we shall weaken this estimate to
$$\Pr[L_{i,j}] = O\left({1\over b}\right).$$
Suppose now that the event $L_{i,j}$ has occurred in some particular way.
Then the displacement from the end of the $k$-th step to the beginning of the $l$-th step
is given by an $a$-step equilateral isotropic random walk, followed by a constant step
of length at most $2$
(from the beginning of the $i$-th step to the end of the $j$-th step), followed by a $c$-step equilateral isotropic random walk.
By the commutativity of addition, this is equivalent to an $(a+c)$-step equilateral isotropic random walk, followed by the constant step.
By the argument in the proof of Proposition 3.1, we obtain
$$\Pr[L_{k,l} \mid L_{i,j}] = {2\over \pi^2(a+c)} + O\left({1\over (a+c)^2}\right),$$
but we shall weaken this estimate to
$$\Pr[L_{k,l} \mid L_{i,j}] = O\left({1\over a+c}\right).$$
Thus
$$\Covar[L_{i,j}, L_{k,l}] = O\left({1\over b}\,{1\over a+c}\right).$$
We shall extend the range of summation to all $1\le a,b,c,d\le n$, which can only increase the 
result.
Summing over $1\le b\le n$ gives a factor of $O(\log n)$, because
$$\sum_{1\le v\le m} {1\over v} = \log m + O(1).$$
The sum over $a$ and $c$ contributes a factor of $O(n)$, because
$$\eqalign{
\sum_{1\le v\le m} \sum_{1\le w\le m} {1\over v+w} 
&= O(m). \cr
}$$
And the sum over $d$ of course contributes a factor of $O(n)$.
Thus the quadruple sum over $a$, $b$, $c$ and $d$ is $O(n^2 \log n)$.
This estimate of course also applies to the sum of terms with $k<l<j<l$.

We consider next the terms with $i<k<j<l$, as depicted below.
\vskip0.15in
$$\xymatrix{
1 \ar@{-}[r] 
&i \ar@{-}[r] _a \ar@{--}@/^2.5pc/[rrr]
&k \ar@{-}[rr] _b \ar@{--}@/_2.5pc/[rrr]
&&j \ar@{-}[r] ^c
&l \ar@{-}[r] ^d
&n \cr
}$$
\vskip0.35in
\centerline{Figure 5.2}
\sk

\noindent For these terms we shall again use the upper bound (5.4), and the variables
$a$, $b$, $c$ and $d$ as defined before.
By Proposition 3.1 
$$\Pr[L_{i,j}] = {2\over \pi^2 (a+b)} + O\left({1\over (a+b)^2}\right),$$
but we shall weaken this estimate to
$$\Pr[L_{i,j}] = O\left({1\over a+b}\right). \eqno(5.5)$$
We claim that
$$\Pr[L_{k,l} \mid L_{i,j}] = O\left({1\over (a\parallel b) + c}\right). \eqno(5.6)$$
To see this, suppose that the event $L_{i,j}$ has occurred in some particular way.
Consider the displacement $A$ from the beginning of the $j$-th step to the end of the
$k$-th step.
This displacement consist of a constant step of distance at most $2$ (from the beginning of the 
$j$-th step to the end of the $i$-th step), followed by an equilateral isotropic random walk of $(a+1)$ steps, conditioned on the event that  a further $b$-step equilateral isotropic random walk  returns to the beginning of the $j$-th step.
It follows that  $A$ is a constant (of length at most $2$) followed by a $(a\parallel b)$-quasi-Gaussian step.
Furthermore, the angular density of the $k$-th step is within a factor $(1 + O(1/(a\parallel b)))$
of uniform.
Consider next the the displacement $B$ from the beginning of the $j$-th step to the beginning of the $l$-th step.
This displacement consists of a constant step of length $1$ (the $j$-th step) followed by a
$c$-step equilateral  isotropic random walk, and is thus a constant step of length $1$ followed by a$c$-quasi-Gaussian step.
Furthermore, the angular density of the $l$-th step is uniform.
We must now consider the total displacement $A+B$ from the end of the $k$-th step to the beginning of the $l$-th step.
Apart from the constant steps, we have a $(a\parallel b)$-quasi-Gaussian step followed by a
$c$-quasi-Gaussian step.
It is well known that the sum of a $v$-Gaussian step and a $w$-Gaussian step is a 
$(v+w)$-Gaussian step.
(The variances of ``series'' Gaussian steps combine like resistance in series.)
For quasi-Gaussian steps, however, we must add error terms $O(1/(v+w)v)$ and $O(1/(v+w)w)$.
These error terms  can be combined to the single error term $O(1/(v+w)(v\parallel w))$.
Thus a sum of quasi-Gaussian steps is not necessarily quasi-Gaussian (that would require
an error term $O(1/(v+w)^2)$), but it differs from quasi-Gaussian only in having the larger error term
$O(1/(v+w)(v\parallel w))$.
Applying this result to the problem at hand, we conclude that
$$\Pr[L_{k,l} \mid L_{i,j}] = {2\over \pi^2((a\parallel b) + c)}
+ O\left({1\over ((a\parallel b) + c)(a\parallel b\parallel c)}\right),$$
but we shall weaken this estimate to (5.6).

Combining (5.4), (5.5) and (5.6), we obtain
$$\eqalign{
\Covar[L_{i,j}, L_{k,l}]  
&= O\left({1\over a+b} \, {1\over (a\parallel b) + c}\right) \cr
& = O\left({1\over ab + ac + bc}\right). \cr
}$$
Since the arithmetic mean $(ab + ac + bc)/3$ exceeds the corresponding geometric mean
$(abc)^{2/3}$, we obtain
$$\Covar[L_{i,j}, L_{k,l}]  = O\left({1\over a^{2/3} \, b^{2/3} \, c^{2/3}}\right).$$
We sum as before over all $1\le a,b,c,d\le n$,
The sums over each of $a$, $b$ and $c$ contribute a factor of $O(n^{1/3})$, because
$$\sum_{1\le v\le n} {1\over v^{2/3}} = O(n^{1/3}).$$
And the sum over $d$ of course contributes a factor of $n$.
Thus the quadruple sum over $a$, $b$, $c$ and $d$ is $O(n^2)$.
This estimate of course also applies to the sum over terms with $k<i<l<j$.

At this point, we have considered all terms in (5.2) in which $i$, $j$, $k$ and $l$ pairwise differ by at least $2$.
The remaining terms are much easier to deal with,
and we shall give the estimates explicitly.
We merely state that these terms 
 contribute only
$O\(n(\log n)^2\)$ to the sum.
Since all of these contributions are $O(n^2 \log n)$, we have verified (5.1).

We turn now to polygons.
Recall that the random variable $K'_n$ denotes the number of self-intersections in a two-dimensional equilateral isotropic random polygon.
We shall derive the estimate
$$\Var[K'_n] = O(n^2 \log n) \eqno(5.7)$$

For polygons, with their circular symmetry, there are only two topologically distinct configurations for a pair of self-intersections, as depicted in Figures 5.4 and 5.5. 
\vskip0.35in
$$\xymatrix{
&&\bullet \ar@{--}[dd]\cr
&\bullet \ar@{-}[ur] \ar@{-}[dr] \ar@{--}@/^5pc/[rr] &&\bullet \ar@{-}[ul] \ar@{-}[dl]\cr
&&\bullet \cr
}$$
\centerline{\ \ \ \ \ \ \ \ \ \ Figure 5.4}
\sk

$$\xymatrix{
&\bullet \ar@{-}[dd] \ar@{-}[rr] \ar@{--}@/^1.5pc/[rr] &&\bullet  \ar@{-}[dd] \cr
\cr
&\bullet \ar@{-}[rr] \ar@{--}@/_1.5pc/[rr] &&\bullet \cr
}$$
\vskip0.35in
\centerline{\ \ \ \ \ \ \ \ \ \ Figure 5.5}
\sk

\noindent
The four stretches  in Figure 5.4 are all topologically equivalent, and drawing any of 
them at the bottom results in a picture like Figure 5.6.
\vskip0.35in
$$\xymatrix{
&i \ar@{-}[r] _a \ar@{--}@/^2.5pc/[rrr] \ar@{-}@/_5pc/[rrrr] ^d
&k \ar@{-}[rr] _b \ar@{--}@/_2.5pc/[rrr]
&&j \ar@{-}[r] ^c
&l  \cr
}$$
\vskip0.25in
\centerline{Figure 5.6}
\sk

\noindent
We shall analyze this case with arguments similar to those we used for Figure 5.3.
But for polygons we have the constraint $a+b+c+d = n-4$, so we shall not sum over
all combinations of $1\le a,b,c,d \le n$.
Rather, we shall take $d = n-4-a-b-c$, sum over all combinations of $1\le a,b,c \le n$,
and then include an extra factor of $n$ to the sum, to take account of the $n$ possible 
positions in which this figure might appear around the polygon.
Furthermore, since at least one of the four stretches in Figure 5.4 must have length
at least $(n-4)/4\ge n/8$ (assuming $n\ge 8$), we may choose such a long stretch 
to draw at the bottom, so that we have $d\ge n/8$.

By the same arguments as we used for Figure 5.3, we have
$$\Pr[L_{i,j}] = O\left({1\over (a+b)\parallel (c+d)}\right)$$
and
$$\Pr[L_{k,l} \mid L_{i,j}] = O\left({1\over (a\parallel b) + (c\parallel d)}\right).$$
Thus
$$\eqalign{
\Covar[L_{i,j}, L_{k,l}] 
&= O\left({1\over (a+b)\parallel (c+d)} \, 
{1\over (a\parallel b) + (c\parallel d)}\right) \cr
\cr
&= O\left({a+b+c+d \over abc+abd+acd+bcd}\right). \cr
}$$
Again using an inequality $(abc+abd+acd+bcd)/4 \ge (abcd)^{3/4}$
between arithmetic and geometric means, together with $a+b+c+d=n-4$, we obtain
$$\Covar[L_{i,j}, L_{k,l}]  = O\left({n\over a^{3/4} \, b^{3/4} \, c^{3/4} \, d^{3/4}}\right).$$
We sum this expression over all $1\le a,b,c\le n$.
The sums over each of $a$, $b$ and $c$ each contribute a factor of $O(n^{1/4})$,
because
$$\sum_{1\le v\le n} {1\over v^{3/4}} = O(n^{1/4}).$$
These contributions are cancelled by the factor of $n^{3/4} \ge (n/8)^{3/4}$
in the denominator.
This leaves just the factor of $n$ in the numerator.
Multiplying by another factor of $n$ to account for the positions in which this configuration
may appear around the polygon, we see that all the terms depicted in Figure 5.5
contribute $O(n^2)$ to the variance.

There are two topologically different kinds of stretches in Figure 5.5.
If we draw it with one of its horizontal stretches at the bottom, we obtain a picture like
Figure 5.7, while if we draw it with one of its vertical stretches at the bottom, we obtain a picture
like Figure 5.8.
\vskip0.35in

$$\xymatrix{
&k \ar@{-}[r] _a \ar@{--}@/^3.75pc/[rrrr] \ar@{-}@/_3.75pc/[rrrr] ^d
&i \ar@{-}[rr] _b \ar@{--}@(ur,ul)[rr]  
&&j \ar@{-}[r] _c
&l \cr
}$$
\vskip0.15in
\centerline{\ \ \ \ \ \ \ \ \ \ \ \ Figure 5.7}
\sk

$$\xymatrix{
&i \ar@{-}[rr] _a \ar@{--}@(ur,ul)[rr] \ar@{-}@/_3.75pc/[rrrrrr] ^d
&&j \ar@{-}[rr] _b
&&k \ar@{-}[rr] _c \ar@{--}@(ur,ul)[rr] 
&&l  \cr
}$$
\vskip0.25in
\centerline{\ \ \ \ \ \ \ \ \ \ \ \ Figure 5.8}
\sk

\noindent
We shall choose between these alternatives to ensure that $d\ge n/8$, as before.

We shall analyze Figure 5.7 with arguments similar to those we used for Figure 5.2.
We shall again take $d = n-4-a-b-c$, sum over all combinations of $1\le a,b,c \le n$,
and then include an extra factor of $n$ to the sum, to take account of the $n$ possible 
positions in which this figure might appear around the polygon.

By the same arguments as we used for Figure 5.2, we have
$$\eqalign{
\Pr[L_{i,j}] 
&= O\left({1\over b\parallel (a+c+d)}\right) \cr
\cr
&= O\left({1\over b}\right) \cr
}$$
and
$$\eqalign{
\Pr[L_{k,l} \mid L_{i,j}] 
&= O\left({1\over (a + c) \parallel  d}\right) \cr
\cr
&= O\left({1\over a + c}\right), \cr
}$$
where we have used the fact that $d\ge n/8$ to simplify the parallel combinations involving $d$.
Thus
$$\Covar[L_{i,j}, L_{k,l}] = O\left({1\over b}\,{1\over a+c}\right).$$
As in the analysis of Figure 5.2, the sum over $1\le b\le n$ contributes a factor of $O(\log n)$
and the sum over $1\le a,c\le n$ contributes a factor of $O(n)$.
Multiplying by another factor of $n$ to account for the positions in which this configuration
may appear around the polygon, we see that all the terms depicted in Figure 5.7
contribute $O(n^2 \log n)$ to the variance.

Finally, we consider Figure 5.8.
This case is most similar to that of Figure 5.1.
But where $L_{i,j}$ and $L_{k,l}$ were independent in Figure 5.1, the stretch of length $d$
introduces a dependence in Figure 5.8.
But since $d\ge n/8$, this dependence is weak.
We shall need to exploit cancellation between the terms in (5.3), rewriting it in the form
$$\Covar[L_{i,j}, L_{k,l}] = \Pr[L_{i,j}] \, \bigg( \Pr[L_{k,l} \mid L_{i,j}]  - \Pr[L_{k,l}] \bigg). \eqno(5.8)$$ 
We have
$$\eqalignno{
\Pr[L_{i,j}]  
&= O\left({1\over a\parallel (b+c+d)}\right) \cr
\cr
&= O\left({1\over a}\right), &(5.9)\cr
}$$

$$\eqalignno{
\Pr[L_{k,l} \mid L_{i,j}]
&= {2\over \pi^2 \(c\parallel (b+d)\)} + O\left({1\over \(c\parallel (b+d)\)^2}\right) \cr
\cr
&= {2\over \pi^2 c} + O\left({1\over c^2}\right) &(5.10)\cr
}$$
and
$$\eqalignno{
\Pr[L_{k,l}]
&= {2\over \pi^2 \(c\parallel (a+b+d)\)} + O\left({1\over \(c\parallel (a+b+d)\)^2}\right) \cr
\cr
&= {2\over \pi^2 c} + O\left({1\over c^2}\right), &(5.11)\cr
}$$
where we have used the fact that $d\ge n/8$ to simplify the parallel combinations involving $d$.
Substituting (5.9), (5.10) and (5.11) in (5.8),
we obtain
$$\Covar[L_{i,j}, L_{k,l}] = O\left({1\over a}\,{1\over c^2}\right).$$
The sum over $1\le a\le n$ contributes a factor of $O(\log n)$, the sum over $1\le b\le n$ contributes a factor of $O(n)$, and the sum over $1\le c\le n$ contributes a factor of $O(1)$.
Multiplying by another factor of $n$ to account for the positions in which this configuration
may appear around the polygon, we see that all the terms depicted in Figure 5.8
contribute $O(n^2 \log n)$ to the variance.

As was the case for walks,  terms in which $i$, $j$, $k$ and $l$ do not differ pairwise by at least $2$ are much easier to deal with,
and contribute only
$O\(n(\log n)^2\)$ to the variance.
Since all of these contributions are $O(n^2 \log n)$, we have verified (5.7).

\heading{6. Conclusion}

We have shown that the mean number of self-intersections for both random walks and random polygons in two dimensions with isotropic equilateral steps is $(2/\pi^2)n\log n + O(n)$.
We have also shown that the variance is $O(n^2 \log n)$ for both walks and polygons.
We have not determined the asymptotic behavior of the variance more exactly, because our result suffices to show concentration around the mean.
It remains an open problem to determine the order of magnitude of (or, more ambitiously, an asymptotic formula for) the variance.

It would also be of interest to extend the results of this paper to Gaussian random walks and polygons, or to the projections onto two dimensions of three-dimensional equilateral walks and polygons,
thereby establishing concentration about the mean for these models.
For these problems, the triple integral evaluated in Section~3 would be replaced by a quintuple integral, and the evaluations of various conditional probabilities would become more complicated,
but the general strategy of our proofs should still be applicable.
\sk

\heading{7. Acknowledgment}

The research reported here was supported
by Grants CCF 0646682 and CCF 0917026 from the National Science Foundation.
\sk

\heading{8. References}

\ref D1; Y. Diao, A. Dobay, R. B. Kusner, K. Millett and A. Stasiak;
``The Average Crossing Number of Equilateral Random Polygons'';
J. Phys.\ A: Math. Gen.; 36 (2003) 11561--11574.

\refinbook D2; Y. Diao and C. Ernst;
``The Average Crossing Number of Gaussian Random Walks and Polygons'';
in J.~A. Calvo, K.~C. Millett and E.~J. Rawdon (editors);
Physical and Numerical Models in Knot Theory;
World Scientific Publishing, 2005, pp.~275--292.

\ref K; Kluyver;
``A Local Probability Problem'';
Proc.\ Section of Sci., K. Akad.\ van Wet., te Amsterdam; 8 (1906) 341--350.

\refbook R1; Lord Rayleigh ($=$ J. W. Strutt);
The Theory of Sound;
Dover Publications, New York, 1945.

\ref R2; Lord Rayleigh ($=$ J. W. Strutt);
``On the Problem of Random Vibrations, and of Random Flights in One, Two, or Three Dimensions'';
Phil.\ Mag.; 37 (1919) 321--347.

\refbook W1; G. N. Watson;
A Treatise on the Theory of Bessel Functions {(\rm second edition)};
Cambridge University Press, New York, 1944.

\refbook W2; E. T. Whittaker and G. N. Watson;
A Course of Modern Analysis {\rm (fourth edition)};
Cambridge University Press, London, 1963.

\bye